\documentclass[12pt]{article}

\usepackage{mathrsfs}
\usepackage{latexsym,amsmath,amssymb,amsfonts,epsfig,cite,psfrag}
\usepackage{eepic,color,colordvi,amscd}
\usepackage{ebezier}
\usepackage{verbatim}
\usepackage{subfigure}
\usepackage{accents}
\usepackage{amsthm}
\usepackage{comment}

\theoremstyle{plain}
\newtheorem{theorem}{Theorem}[section]
\newtheorem{lemma}[theorem]{Lemma}
\newtheorem{corollary}[theorem]{Corollary}

\newcommand{\claimqed}{{\hfill \rule{4pt}{7pt}}}

\newcommand{\ceil}[1]{\left\lceil#1\right\rceil}

\newcommand{\define}[1]{\emph{#1}}

\newcommand{\bridges}[2][]{\beta_{#1}(#2)}
\newcommand{\taubridges}[2][]{\tau_{#1 #2}}

\newcommand{\bridgecorrection}[2][]{\phi_{#1}(#2)}

\newcommand \track[1] {{\color{blue} #1}}
\newcommand \sout[1] {{\color{red} #1}}
\renewcommand\track[1]{{\color{black} #1}}
\renewcommand\sout[1]{{\iffalse #1 \fi}}
\begin{document}

\title{On Tutte cycles containing three prescribed edges}
\author{Michael C. Wigal\footnote{Department of Mathematics, University of Illinois at Urbana--Champaign, Urbana, IL, USA. Supported by an NSF Graduate
    Research Fellowship under Grant No. DGE-1650044. Email: wigal@illinois.edu}   \ \ and Xingxing Yu\footnote{	School of Mathematics,
    Georgia Institute of Technology,
    Atlanta, GA, USA. Partially
    supported by NSF Grant DMS 1954134. Email: yu@math.gatech.edu}\\
}
\date{\today}

\maketitle

\begin{abstract}
	
	A cycle $C$ in a graph $G$ is called a Tutte cycle if, after deleting $C$ from $G$, each component has at most three neighbors on $C$. Tutte cycles play an important role in the study of Hamiltonicity of planar graphs. Thomas and Yu and independently Sanders proved the existence of Tutte cycles containining three specified edges of a facial cycle in a 2-connected plane graph. We prove a quantitative version of this result, bounding the number of components of the graph obtained by deleting a Tutte cycle. As a corollary, we can find long cycles in essentially 4-connected plane graphs that  also contain three prescribed edges of a facial cycle. \\
	
\end{abstract}

\newpage

\section{Introduction}\label{section_1}

Let $G$ be a graph. We use $V(G)$ and $E(G)$ to denote the vertex set and edge set of $G$, respectively, and use $|G|$ to denote the number of vertices in $G$. Let $H$ be a subgraph of $G$. We use $G-H$ to denote the graph obtained from $G$ by deleting $H$ and all edges of $G$ incident with $H$.  An $H$-\define{bridge} of $G$ is a subgraph of $G$ that is either induced by an edge in $E(G) \setminus E(H)$ with both incident vertices belonging to $H$, or induced by the edges of $G$ incident with at least one vertex of a single component of $G - H$. (An $H$-bridge induced by a single edge \track{with both incident vertices belonging to $H$} is \define{trivial}.) Given an $H$-bridge $B$ of $G$, the \define{attachments} of $B$ on $H$ are \track{the} vertices in $V(H) \cap V(B)$.  A subgraph $H$ of $G$ is \define{Tutte} if every $H$-bridge of $G$ has at most three attachments on $H$. If $G$ is 4-connected with $|G| > 4$ and $H$ is a Tutte subgraph of $G$, then $V(G) = V(H)$.

Tutte subgraphs are important in the study of Hamiltonicity of planar graphs. 
A path or cycle in a graph $G$ is \define{Hamiltonian} if it visits every vertex of $G$. A graph is \define{Hamiltonian} if it contains a Hamiltonian cycle, and is \define{Hamiltonian-connected} if there is a Hamiltonian path connecting every pair of vertices. 
Tait \cite{T84} conjectured that 3-connected \track{cubic} planar graphs are Hamiltonian. Tutte \cite{T46} disproved Tait's conjecture with a counterexample, and there has been work on finding more, see \cite{ABHM00, HM88}.  Later, Tutte \cite{T56} showed the existence of Tutte cycles in 2-connected planar graphs. As a corollary, all 4-connected planar graphs are Hamiltonian.

Let $H$ and $S$ be subgraphs of a graph $G$. $H$ is said to be an $S$-\define{Tutte subgraph} if $H$ is Tutte and every $H$-bridge of $G$ containing an edge of $S$ has at most two attachments on $H$. \sout{Thomassen \cite{T83} proved that if $G$ is a 2-connected plane graph and $C$ is a facial cycle of $G$ then $G$ has a $C$-Tutte path between two given vertices of $C$ and containing a specified edge of $C$.} 

\track{Let $G$ be a 2-connected plane graph and $C$ be the outer cycle of $G$, i.e., the cycle bounding the infinite face of $G$. Let $v \in V(C)$, $e \in E(C)$, and $u \in V(G)$ such that $u$ and $v$ are distinct. Thomassen \cite{T83}, with a small correction by Chiba and Nishizeki \cite{CN86}, proved that there exists an $C$-Tutte path $P$ between $u$ and $v$ such that $e \in E(P)$.} As a consequence, 4-connected planar graphs are Hamiltonian-connected.

Wigal and Yu 
\cite{WY20} proved a quantitative version of Thomassen's theorem by controlling the number of bridges. Our goal is to extend this work to graphs on other surfaces, in particular, the projective plane. As a first step towards this goal we prove a quantitative version of a result proved independently by \sout{Sander}\track{Sanders} \cite{S96} and Thomas and Yu \cite{TY94}. 
In order to state our result, we need some notation. For a positive integer $k$, \sout{A}\track{a} $k$-\define{separation} in a graph $G$ is a pair of subgraphs of $G$, say $(G_1,G_2)$, such that $|V(G_1) \cap V(G_2)| = k$, $E(G_1) \cap E(G_2) = \emptyset$, $G_1 \cup G_2 = G$, and $G_i \not \subseteq G_{3 - i}$ for $i \in \{1,2\}$. \track{Note that we allow for the possibility of either $G_1$ or $G_2$ consisting of a single edge.} A $k$-\define{cut} in $G$ is a set $S \subseteq V(G)$ such that $|S| = k$ and there exists a separation $(G_1,G_2)$ such that $V(G_1) \cap V(G_2) = S$ and \sout{$G_i - G_{3-i} \neq \emptyset$}\track{$V(G_i) \setminus V(G_{3-i}) \neq \emptyset$} for $i \in \{1,2\}$. 

The \define{outer walk} of a plane graph $G$ consists of vertices and edges of $G$ incident with its infinite face, which is called \track{the} \define{outer cycle} if it is a cycle. A \define{circuit graph} $(G,C)$ consists of a 2-connected plane graph $G$ and outer cycle $C$, such that for any 2-cut $T$ of $G$, each component of $G - T$ must contain an edge of $E(C)$. Let $G$ be a 2-connected plane graph with \sout{facial}\track{outer} cycle $C$. For $x,y \in V(C) \cup E(C)$, we define $xCy$ as the subpath of $C$ in clockwise order from $x$ to $y$ such that $x,y \not \in E(xCy)$. If $x = y \in V(C)$, then $xCy$ consists of the vertex $x = y$. We say that $xCy$ is \define{good} if $G$ has no 2-separation $(G_1,G_2)$ with $V(G_1 \cap G_2) = \{s,t\}$ such that $sCt \subseteq xCy \cap G_2$ and $|G_2| \ge 3$. \track{Let $H \subseteq G$ and let
	\begin{align*}
		\bridges[G]{H} = |\{ B : \text{$B$ is an $H$-bridge and $|B| \ge 3$}\}|.
\end{align*}} For $x,y \in V(C) \cup E(C)$, define
\[
\taubridges[G]{xy} =
\begin{cases}
	2/3, \quad & \text{$xCy$ is not good,}\\
	2/3 \quad & \text{$|\{x,y\}\cap E(C)| \ge 1$ and $x$ and $y$ are incident, }\\
	1/3, \quad & \text{$|\{x,y\}\cap E(C)| \ge 1$ and $|xCy|=2$,}\\
	0, \quad & \text{otherwise,}
\end{cases}
\]
and, for any $P \subseteq G$, define
\[\bridgecorrection[G]{P} =  \sum_{B \in {\cal B}_P}\frac{(|B| - 3)}{3},\]
where ${\cal B}_P$ is the set of all nontrivial $P$-bridges of $G$ with 2 attachments on $P$. When there is no ambiguity, we drop the reference to $G$ and simply write \track{$\bridges{H}$,} $\taubridges{xy}$, or $\bridgecorrection{P}$. Our main result is a quantitative version of a result by Sanders \cite{S96} and Thomas and Yu \cite{TY94}, on Tutte cycles containing three prescribed edges on a facial cycle. 
\begin{theorem}\label{two_edge_tech}
	Let $n \ge 3$ be an integer, let $(G,C)$ be a circuit graph on $n$ vertices, let $u,v \in V(C)$ be distinct with $uv \in E(C)$, and let $e ,f \in E(C)$ be distinct, such that $u,f,e,v$ occur on $C$ in clockwise order. Then $G$ has a $C$-Tutte path $P$ between $u$ and $v$, such that $\{e,f \}\subseteq  E(P)$, and $\bridges{P} \le (n - 7)/3 + \taubridges{uf} + \taubridges{fe}  + \taubridges{ev} - \bridgecorrection{P}$.
\end{theorem}

In Section \ref{section_2}, we state the main result from  \cite{WY20} and prove additional lemmas.  We prove Theorem \ref{two_edge_tech} in Section \ref{section_3}.  In Section \ref{section_4}, we discuss applications of Theorem \ref{two_edge_tech} to finding long cycles in planar graphs that are not 4-connected and some immediate corollaries. \sout{We also prove a corollary  of Theorem \ref{two_edge_tech} that will be of use to extend our work to projective-planar graphs.}

We conclude this section with additional notation. Let $G$ be a graph. If $T$ is a collection of 2-element subsets of $V(G)$, we let $G + T$ denote the graph on $V(G)$ with edge set $E(G) \cup T$. If $T = \{\{u,v\}\}$, we write $G + uv$ \sout{as well}\track{for $G + T$}. It is \sout{convient}\track{convenient} to regard paths as a sequence of vertices, with \sout{consectutive}\track{consecutive} vertices adjacent. If $H$ is also a graph, we let $G \cup H$ and $G \cap H$ denote the union and intersection of $G$ and $H$, respectively.
\section{Preliminary results}\label{section_2}

We begin with the following result of Wigal and Yu \cite{WY20}, which is a quantitative version of a result of Thomassen \cite{T83}.

\begin{lemma}\label{planar_technical}
	Let $n\ge 3$ be an integer, let $(G,C)$ be a circuit graph on 
	$n$ vertices, and let $u,v \in V(C)$ be distinct, and
	let $e\in  E(C)$. Suppose that $u,e, v$ occur on $C$ in clockwise order.
	Then $G$ has a $C$-Tutte path $P$ between $u$ and $v$, such that $e\in E(P)$ and
	\[\bridges{P} \le (n - 6)/3 + \taubridges{vu} + \taubridges{ue} + \taubridges{ev}\]
\end{lemma}

From the proof of Lemma \ref{planar_technical} in \cite{WY20}, we see that it will be convenient for counting purposes to assign weight to bridges with two attachments: For each bridge $B$ with two \sout{attachment}\track{attachments}, we wish to only count one internal vertex; this prevents over counting $(|B| - 3)/3$ vertices in the induction step. In the multiple places where we apply Lemma \ref{planar_technical}, it is convenient to work directly with the following ``refined" version using $\bridgecorrection{P}$.

\begin{lemma}\label{planar_technical_refined}
	Let $n \ge 3$ be an integer, let $(G,C)$ be a circuit graph on 
	$n$ vertices, let $u,v \in V(C)$ be distinct, and
	let $e \in  E(C)$. Suppose $u,e, v$ occur on $C$ in clockwise order.
	Then $G$ has a $C$-Tutte path $P$ between $u$ and $v$, such that $e\in E(P)$ and
	$\bridges{P} \le (n - 6)/3 + \taubridges{vu} + \taubridges{ue} + \taubridges{ev} -\bridgecorrection{P}$.
\end{lemma}

\begin{proof}
	A 2-separation $(G_1,G_2)$ in $G$ with $\{u,v\} \subseteq V(G_1)$ and $e \in E(G_1)$ is said to be
	\define{maximal} if there is no 2-separation $(G_1',G_2')$ in $G$ such that  $\{u,v\} \subseteq V(G_1')$, $e \in E(G_1')$, $G_2 \subseteq G_2'$, and $G_2 \neq G_2'$.
	
	For every maximal 2-separation $(G_1,G_2)$ of $G$,  we contract \sout{$G_2$}\track{$G_2 - V(G_1)$} to a new vertex, so that all these new vertices are distinct. Let $T$ be the set of new vertices. Let $H$ denote the resulting 2-connected plane graph, and let $D$ denote the outer cycle of $H$. Note that $v,u \in V(H)$ and $e \in E(H)$. Furthermore, it is straightforward to check that  $\taubridges[G]{vu} = \taubridges[H]{vu}$, $\taubridges[G]{ue} = \taubridges[H]{ue}$, and  $\taubridges[G]{ev} = \taubridges[H]{ev}$. 
	Applying Lemma \ref{planar_technical} to $H$, we obtain a $D$-Tutte path $P'$ in $H$ between $u$ and $v$, such that $e \in E(P)$ and
	\[\bridges[H]{P'} \le (|H| - 6)/3 + \taubridges[H]{vu} + \taubridges[H]{ue} + \taubridges[H]{ev} =  (|H| - 6)/3 + \taubridges[G]{vu} + \taubridges[G]{ue} + \taubridges[G]{ev}.\]
	To extend $P'$ to the desired path $P$ in $G$, we need to handle vertices in $T \cap V(P')$. Let $t \in T \cap V(P')$ with neighbors $x$ and $y$. Then there exists a 2-separation $(G_1,B_t)$ in $G$ with $v,u \in V(G_1)$, $e \in E(G_1)$, and $V(G_1 \cap B_t) = \{x,y\}$. Without loss of generality, assume $xCy \subseteq B_t$. Let $B_t' = B_t + xy$ and $C_t = xCy + yx$. Then we may draw $B_t'$ in the plane so that $(B_t',C_t)$ is a circuit graph.  Since $(G,C)$ is a circuit graph, $|xCy|\ge 3$. Choose an edge $e' \in xCy$ such that \sout{$\taubridges[B_t']{xe} \le 1/3$}\track{$\taubridges[B_t']{xe'} \le 1/3$, e.g., take $e'$ to be the second edge of $xCy$.} Note that $\taubridges[B'_t]{yx}=0$. Then by induction, there exists a $C_t$-Tutte path \sout{$P_t$ in $B'_t$}\track{$P_t \subseteq B'_t$ between $x$ and $y$} such that \track{$e' \in P_t$ and} \begin{align*}
		\bridges[B'_t]{P_t} &\le (|B'_t| - 6)/3 + \taubridges[B'_t]{yx} + \taubridges[B'_t]{xe'} + \taubridges[B'_t]{e'y} -\bridgecorrection[B'_t]{P_t}\\
		&\le (|B'_t| - 3)/3 - \bridgecorrection[B'_t]{P_t}\\
		&= (|B_t| - 3)/3 - \bridgecorrection[B_t]{P_t}.
	\end{align*}
	Now $P = (P' - T) \cup (\bigcup_{t \in T \cap V(P')} P_t)$ is a $C$-Tutte path in $G$ between $u$ and $v$ such that $e \in E(P)$. Next, we bound $\bridges[G]{P}$. For every $t \in T \backslash V(P')$, there is a corresponding $P'$-bridge $B_t$ of $G$ with $|B_t| - 3$ vertices of $G$ not counted in $|H|$, where we count $t$ as one vertex from $B_t$. For each $t \in T \cap V(P')$, $|B_t|$ and $|H|$ overcount by three vertices, the two vertices in $V(B_t \cap P)$ and the vertex $t$. Hence,
	\begin{align*}
		\bridges[G]{P} &\le (|H| - 6)/3 + \taubridges[G]{vu} + \taubridges[G]{ue} + \taubridges[G]{ev}\\
		&+\sum_{t \in V(P') \cap T}((|B_t| - 3)/3 - \bridgecorrection[B_t]{P_t})\\
		&= (|G| - 6)/3 + \taubridges[G]{vu} + \taubridges[G]{ue} + \taubridges[G]{ev}\\
		&- \sum_{t \in V(P') \cap T} \bridgecorrection[B_t]{P_t} - \sum_{t \in T \setminus V(P')} (|B_t| - 3)/3\\
		&= (n - 6)/3 + \taubridges[G]{vu} + \taubridges[G]{ue} + \taubridges[G]{ev} -\bridgecorrection[G]{P}
	\end{align*} 
	as $\bridgecorrection[G]{P} = \sum_{t \in V(P') \cap T} \bridgecorrection[B_t]{P_t} + \sum_{t \in T \setminus V(P')} (|B_t| - 3)/3$. 
\end{proof}

Next, we derive a corollary of Lemma \ref{planar_technical_refined}.

\begin{corollary}\label{planar_technical_vertex}
	Let $n \ge 3$ be an integer, let $(G,C)$ be a circuit graph on $n$ vertices, and let $u,z,v \in V(C)$ be distinct \track{and occur on $C$ in clockwise order}. Then $G$ has a $C$-Tutte path $P$ between $u$ and $v$, such that $z \in V(P)$ and $\bridges{P} \le (n - 3)/3 + \taubridges{vu} -\bridgecorrection{P}$.
\end{corollary}

\begin{proof}
	\sout{Without loss of generality, assume that $u,z,v$ occur on $C$ in clockwise order.} If there is not a 2-cut in $G$ separating $z$ from $\{v,u\}$, then we choose an edge $e$ on $uCv$ such that $\taubridges{ue} \le 1/3$ and $\taubridges{ev} \le 2/3$. Applying Lemma \ref{planar_technical_refined} we find a $C$-Tutte path $P$ in $G$ between $u$ and $v$, such that $e \in E(P)$ and
	\[ \bridges[G]{P} \le (n - 3)/3 + \taubridges[G]{vu} -\bridgecorrection[G]{P}.\]
	Note $z \in V(P)$, since there is no 2-cut in $G$ separating $z$ from $v$ and $u$. So $P$ is the desired path.
	
	Now assume that there exists a 2-separation $(G_1,G_2)$ in $G$ with $\{u,v\} \subseteq V(G_1)$ and $z \in V(G_2)$. Let $V(G_1 \cap G_2) = \{x,y\}$ such that $u,x,z,y,v$ appear on $C$ in clockwise order. Furthermore, let $G_1' := G_1 + xy$ and $G_2' := G_2 + xy$ be plane graphs such that $G_1'$ has outer cycle $C'_1 := yCx + xy$ and $G_2'$ has outer cycle $C'_2 := xCy + xy$. Then $(G_i',C_i')$, $i = 1,2$, are circuit graphs.
	
	By Lemma \ref{planar_technical_refined}, $G_1'$ has a $C'_1$-Tutte path $P_1$ between $u$ and $v$, such that  $e' := xy \in E(P_1)$ and
	\begin{align*}
		\bridges[G_1']{P_1} &\le (|G_1'| - 6)/3 + \taubridges[G_1']{vu} + \taubridges[G_1']{ue'} + \taubridges[G_1']{e'v} - \bridgecorrection[G_1']{P_1}\\
		&\le (|G_1'| -2)/3 + \taubridges[G]{vu} - \bridgecorrection[G_1']{P_1}.
	\end{align*}
	By induction, $G_2'$ has a $C'_2$-Tutte path $P_2$ between $x$ and $y$ such that $z \in V(P)$ and 
	\[ \bridges[G_2']{P_2} \le (|G_2'| -3)/3 +\taubridges[G_2']{yx}- \bridgecorrection[G_2']{P_2} = (|G_2'| -3)/3 - \bridgecorrection[G_2']{P_2}.\]
	Let $P = (P_1 - xy) \cup P_2$. Since $|G_1'| + |G_2'| = n + 2$ and $\bridgecorrection[G]{P} = \bridgecorrection[G_1']{P} + \bridgecorrection[G_2']{P_2}$, we see that $P$ is a $C$-Tutte path between $u$ and $v$ in $G$, such that $z \in V(P)$ and $ \bridges[G]{P} \le (n - 3)/3 + \taubridges[G]{vu} -\bridgecorrection[G]{P}$. 
\end{proof}

We now prove a lemma that will help us deal with special 2-separations in the proof of Theorem \ref{two_edge_tech}.

\begin{lemma}\label{two_edge_cut}
	Suppose $n \ge 4$ is an integer and Theorem \ref{two_edge_tech} holds for graphs on at most $n - 1$ vertices. Let $(G,C)$ be a circuit graph on $n$ vertices, let $u,v \in V(C)$ be distinct with $uv \in E(C)$, and let $e,f \in E(C)$ be distinct, such that $u,f,e,v$ occur on $C$ in clockwise order.
	
	Suppose $G$ has a 2-separation $(G_1,G_2)$ such that $\{u,v\} \subseteq V(G_1), \{u,v\} \neq V(G_1 \cap G_2)$, $e \in E(G_2)$, and $|G_2| \ge 3$. Then $G$ has a $C$-Tutte path $P$ between $u$ and $v$, such that $e,f \in E(P)$ and $\bridges[G]{P} \le (n - 7)/3  + \taubridges[G]{uf}+ \taubridges[G]{fe} + \taubridges[G]{ev} - \bridgecorrection[G]{P}$.
\end{lemma}

\begin{proof}
	Let $V(G_1 \cap G_2) = \{x,y\}$ and assume $v,u,x,e,y$ occur on $C$ in clockwise order. Let $G'_i = G_i + xy$ for $i \in \{1,2\}$, such that $G_1'$ is a
	plane graph with outer cycle $C_1:=yCx + xy$ and $G_2'$ is a plane graph with
	outer cycle $C_2:=xCy+xy$. Let
	$e_1:=xy$. Note $|G_1'| + |G_2'| = n + 2$, and both $(G_1',C_1)$ and $(G_2',C_2)$ are circuit graphs. Since $uv \in E(C_1)$,  $\taubridges[G_1']{vu} = 0$.\\
	
	\textit{Case} 1. $f \in E(G_2)$.
	
	By Lemma \ref{planar_technical_refined}, $G_1'$ has a $C_1$-Tutte path $P_1$  between $u$ and $v$,  such that $e_1 \in E(P_1)$ and 
	\[ \bridges[G'_1]{P_1} \le (|G_1'| - 6)/3  + \taubridges[G_1']{ue_1}  + \taubridges[G_1']{e_1v} -\bridgecorrection[G_1']{P_1}.\]
	By assumption, $G_2'$ has a $C_2$-Tutte path $P_2$ between $x$ and $y$, such that $\{e,f\} \subseteq E(P)$ and
	\[\bridges[G_2']{P_2} \le (|G_2'|-7)/3 + \taubridges[G_2']{xf} + \taubridges[G_2']{fe} + \taubridges[G_2']{ey} -\bridgecorrection[G_2']{P_2}.\]
	Note that $P = (P_1 - e_1) \cup P_2$ is a $C$-Tutte path in $G$ between $u$ and $v$, $\{e,f\}\subseteq E(P)$, and $\bridgecorrection[G_1']{P_1} + \bridgecorrection[G_2']{P_2} = \bridgecorrection[G]{P}$. Moreover, 
	\begin{align*}
		\bridges[G]{P} &= \bridges[G'_1]{P_1} + \bridges[G_2']{P_2}\\
		&\le  (|G_1'| - 6)/3 + \taubridges[G_1']{ue_1}+ \taubridges[G_1']{e_1v}  -\bridgecorrection[G_1']{P_1}\\
		&+  (|G_2'|-7)/3 + \taubridges[G_2']{xf} + \taubridges[G_2']{fe} + \taubridges[G_2']{ey} -\bridgecorrection[G_2']{P_2}\\
		&= (n - 7)/3 - 4/3 +\taubridges[G_1']{ue_1} + \taubridges[G_1']{e_1v}  + \taubridges[G_2']{xf} + \taubridges[G_2']{fe} + \taubridges[G_2']{ey} -\bridgecorrection[G]{P}.
	\end{align*}
	We claim that $\taubridges[G_1']{ue_1} + \taubridges[G_2']{xf} \le 2/3 + \taubridges[G]{uf}$ and \sout{$\taubridges[G_1']{e_1v} + \taubridges[G_2']{ey} \le 2/3 + \taubridges{ev}$} \track{$\taubridges[G_1']{e_1v} + \taubridges[G_2']{ey} \le 2/3 + \taubridges[G]{ev}$.} By symmetry we only prove the former. If $u = x$, then $\taubridges[G_1']{ue_1} = 2/3$ and $\taubridges[G_2']{xf} = \taubridges[G]{uf}$ thus the inequality holds; if $uCf$ is not good, then $\taubridges[G]{uf} = 2/3$, so the inequality holds. Thus assume $u\ne x$ and $uCf$ is good. If $\taubridges[G_1']{ue_1} = 0$, again the inequality holds; thus assume $\taubridges[G_1']{ue_1} = 1/3$. If $\taubridges[G_2']{xf} \le 1/3$, the inequality holds; thus assume $\taubridges[G_2']{xf} = 2/3$. Since $uCf$ is good, $x$ is incident with $f$. This implies $\taubridges[G]{uf} = \taubridges[G_1']{ue_1} = 1/3$; again the inequality holds. 
	
	Therefore, since $\taubridges[G_2']{fe} = \taubridges[G]{fe}$, we have
	\begin{align*}
		\bridges[G]{P} &\le (n - 7)/3 - 4/3 + (2/3 + \taubridges{uf}) + (2/3 + \taubridges{ev})  + \taubridges{fe}- \bridgecorrection[G]{P}\\
		&= (n - 7)/3  + \taubridges{uf}  + \taubridges{fe}+ \taubridges{ev}- \bridgecorrection[G]{P}.
	\end{align*}
	
	\textit{Case} 2. $f \in E(G_1)$. 
	
	By assumption, $G_1'$ has a $C_1$-Tutte path $P_1$ between $u$ and $v$, such that $\{e_1,f\}\subseteq  E(P_1)$ and
	\[ \bridges[G_1']{P_1} \le (|G_1'| - 7)/3 +\taubridges[G_1']{uf} + \taubridges[G_1']{fe_1} + \taubridges[G_1']{e_1v} - \bridgecorrection[G_1']{P_1}.\]
	By Lemma \ref{planar_technical_refined}, $G_2'$ has an $C_2$-Tutte path $P_2$ between $x$ and $y$, such that $ e \in E(P)$ and
	\[ \bridges[G_2']{P_2} \le (|G_2| -6)/3 + \taubridges[G_2']{yx} + \taubridges[G_2']{xe} + \taubridges[G_2']{ey} - \bridgecorrection[G_2']{P_2}.\]
	Now $P := (P_1 - e_1) \cup P_2$ is a $C$-Tutte path in $G$ between $u$ and $v$, such that $\{e,f\}\subseteq  E(P)$. Moreover,  since $\taubridges[G_2']{yx} = 0$ and $\taubridges[G_1']{uf} = \taubridges[G]{uf}$, we have
	\begin{align*}
		\bridges[G]{P} &=\bridges[G_1']{P_1} + \bridges[G_2']{P_2}\\
		&\le  (|G_1'| - 7)/3 +\taubridges[G_1']{uf} + \taubridges[G_1']{fe_1} + \taubridges[G_1']{e_1v} - \bridgecorrection[G_1']{P_1}\\
		&+ (|G_2'| -6)/3 + \taubridges[G_2']{yx} + \taubridges[G_2']{xe} + \taubridges[G_2']{ey} - \bridgecorrection[G_2']{P_2}\\
		&= (n - 7)/3 -4/3 + \taubridges[G]{uf} + \taubridges[G_1']{fe_1} + \taubridges[G_1']{e_1v} + \taubridges[G_2']{xe} + \taubridges[G_2']{ey} - \bridgecorrection[G]{P}.
	\end{align*}
	
	\sout{Since $\taubridges[G_1']{fe_1} = \taubridges[G_1']{fx}$, as in Case 1 we can show that that $\taubridges[G_1']{fe_1} + \taubridges[G_2']{xe} \le 2/3 + \taubridges[G]{fe}$. In a similar manner as $\taubridges[G_1']{e_1v} = \taubridges[G_1']{yv}$, we can show that $\taubridges[G_1']{e_1v} + \taubridges[G_2']{ey} \le 2/3 + \taubridges[G]{ev}$.}
	\track{By the same argument as in Case 1, we have $\taubridges[G_1']{e_1v} + \taubridges[G_2']{ey} \le 2/3 + \taubridges[G]{ev}$. Now let $w \in fCe$ be the vertex incident with $f$. Again, by the same argument as in Case 1, we have $\taubridges[G_1']{we_1} + \taubridges[G_2']{xe} \le 2/3 + \taubridges[G]{we}$. As $\taubridges[G_1']{we_1} = \taubridges[G_1']{fe_1}$ and $\taubridges[G]{we} = \taubridges[G]{fe}$, we have, $\taubridges[G_1']{fe_1} + \taubridges[G_2']{xe} \le 2/3 + \taubridges[G]{fe}$.} Hence,
	\begin{align*}
		\bridges[G]{P} &\le (n - 7)/3 -4/3+ \taubridges[G]{uf}  +(2/3 + \taubridges[G]{fe}) + (2/3 + \taubridges[G]{ev}) - \bridgecorrection[G]{P}\\
		&=  (n - 7)/3 + \taubridges[G]{uf} + \taubridges[G]{fe}+ \taubridges[G]{ev}  - \bridgecorrection[G]{P}.
	\end{align*}
	\sout{Hence,}\track{Thus we may conclude} $P$ gives the desired path. 
\end{proof}

\section{Proof of Theorem \ref{two_edge_tech}}\label{section_3}

We apply induction on $n$. The assertion clearly holds for $n = 3$. Hence, we may assume $n > 3$ and the assertion of Theorem \ref{two_edge_tech} holds for graphs on fewer than $n$ vertices. We may assume that

\begin{itemize}
	\item[(1)] $G$ has no  2-separation $(G_1,G_2)$ such that $\{u,v\} \subseteq V(G_1)$, $\{e,f\} \subseteq E(G_1)$, and $|G_2| > 3$. 
\end{itemize}
For otherwise, let $(G_1,G_2)$ be such a 2-seperation in $G$. Let $V(G_1 \cap G_2) = \{x,y\}$ such that $G_2 \cap C = xCy$. Let $G'$ be obtained from $G - (G_2  -\{x,y\})$ by adding a new vertex $t$ and two edges $tx,ty$, such that $G'$ is a plane graph with outer cycle $C' := yCx \cup xty$. Then $(G',C')$ is a circuit graph.  Applying induction to $(G',C')$, we find a $C'$-Tutte path $P'$ in $G'$ between $u$ and $v$, such that $\{e,f\} \subseteq E(P')$ and
\[ \bridges[G']{P'} \le (|G'| - 7)/3 + \taubridges[G']{uf} + \taubridges[G']{fe} + \taubridges[G']{ev} - \bridgecorrection[G']{P'}.\] 
Note that $\taubridges[G']{uf} = \taubridges[G]{uf}$, $\taubridges[G']{fe} = \taubridges[G]{fe}$, and $\taubridges[G']{ev} = \taubridges[G]{ev}$. If $t \not \in V(P')$ then $P = P'$ is also a $C$-Tutte path in $G$.\sout{ and $|G|/3 - \bridgecorrection[G]{P} = |G'|/3 - \bridgecorrection[G']{P'}$; thus $P$ is the desired path. } \track{Furthermore, we have
	\[ \bridgecorrection[G]{P} - \bridgecorrection[G']{P'} = (|G_2| - 3)/3 = (|G| -|G'|)/3.\] Thus we have 
	\[ \bridges[G]{P} \le (n-7)/3 + \taubridges[G]{uf} + \taubridges[G]{fe}+ \taubridges[G]{ev}  - \bridgecorrection[G]{P}\] 
	and conclude $P$ is the desired path.}

Now assume $t \in V(P')$. We view $\track{G_2' :=}G_2+xy$ \track{as a} plane graph with outer cycle $C_2:=xCy+yx$. Then $(\sout{G_2 + xy}\track{G_2'},C_2)$ is a circuit graph. By Corollary \ref{planar_technical_vertex}, $G_2'$ has a $C_2$-Tutte path $P_2$ between $x$ and $y$, such that
\[ \bridges[\sout{G_2}\track{G_2'}]{P_2} \le (|\sout{G_2}\track{G_2'}| - 3)/3 - \bridgecorrection[\sout{G_2}\track{G_2'}]{\sout{P_t}\track{P_2}}. \]
Note that $|G'| + |\sout{G_2}\track{G_2'}| = n + 3$ and that $P = (P' - t) \cup P_2$ is a $C$-Tutte path in $G$ between $u$ and $v$ such that $\{e,f\} \subseteq E(P)$. Moreover, $\bridgecorrection[\sout{G_2}\track{G_2'}]{\sout{P_t}P_2}+\bridgecorrection[G']{P'}=\bridgecorrection[G]{P}$. Hence, 
\begin{align*}
	\bridges[G]{P} &= \bridges[G']{P'} + \bridges[\sout{G_2}\track{G_2'}]{\sout{P_t}\track{P_2}}\\
	&\le (|G'| - 7)/3 + \taubridges[G']{uf} + \taubridges[G']{fe} + \taubridges[G']{ev} - \bridgecorrection[G']{P'}\\
	&+  (|\sout{G_2}\track{G_2'}| - 3)/3 - \bridgecorrection[\sout{G_2}\track{G_2'}]{\sout{P_t}\track{P_2}}\\
	&= (n - 7)/3 + \taubridges[G]{uf} + \taubridges[G]{fe} + \taubridges[G]{ev} -\bridgecorrection[G]{P}.
\end{align*}
Thus $P$ is the desired path; so we may assume (1). 
\claimqed
\medskip

\track{From (1) we see that $\bridgecorrection{P} = 0$ for any Tutte path $P$ in $G$ from $u$ to $v$ trough $e$ and $f$.}

\begin{itemize}
	\item[(2)] We may assume that $|uCe| \ge 3$ or $|fCv| \ge 3$. 
\end{itemize}
For, otherwise, $|C|=3$ and $\taubridges{uf}= \taubridges{fe} = \taubridges{ev}=2/3$. So $P = C-uv$ is a $C$-Tutte path in $G$ through $e$ and $f$\track{.}\sout{ and $\bridgecorrection[G]{P}=1$. Since $n \ge 4$, we have} \track{As $n \ge 4$ and $|C| = 3$, we have that $\bridges[G]{P} = 1$. It follows,}
\begin{align*}
	\bridges[G]{P} = 1 \le (n-7)/3 + \taubridges{uf} + \taubridges{fe} + \taubridges{ev} - \bridgecorrection{P}.
\end{align*}
Thus,  $P$ is the desired path. 

\claimqed

\medskip

\begin{center}
\begin{figure}
	\hspace*{2.5cm} \includegraphics[scale=0.35]{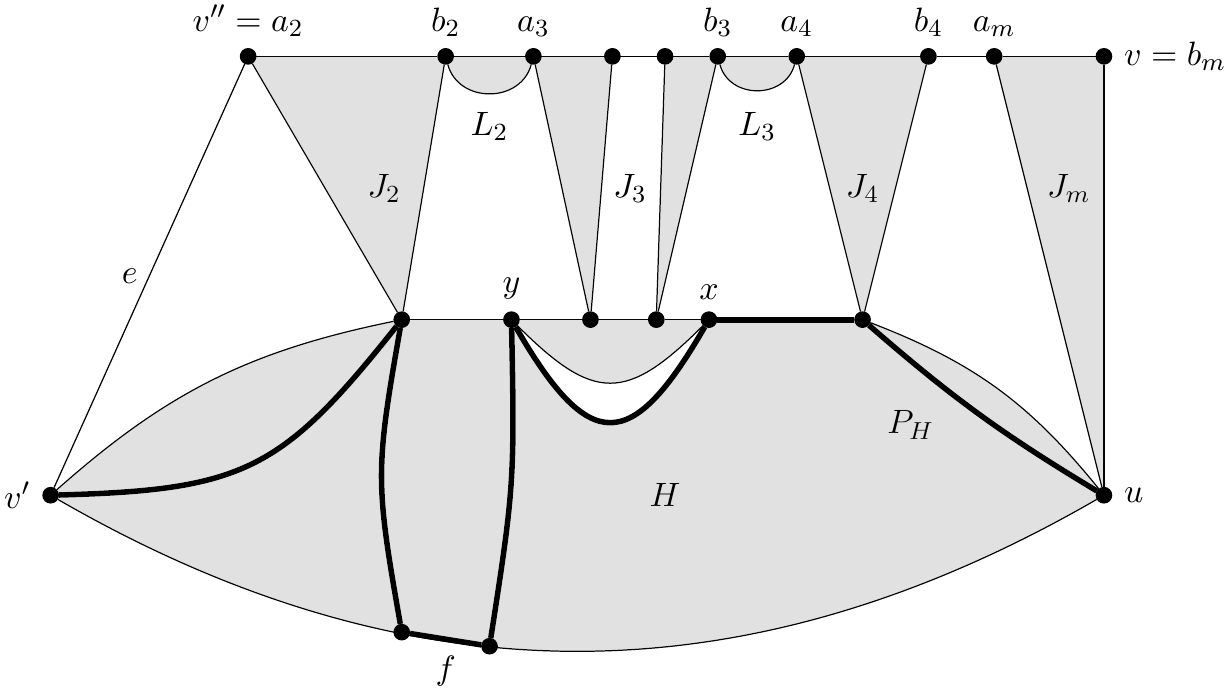}
	\caption{The subgraph $H$, path $P_H$, and some  $(eCv \cup H)$-bridges of $G$.}
	\label{fig:1}
\end{figure}
\end{center}

\sout{From (1) we see that $\bridgecorrection{P} = 0$ for any Tutte path $P$ in $G$.} By (2) and symmetry, we may assume $|uCe| \ge 3$.  Let $e = v'v''$ such that $u,f,v',v'',v$ appear on $C$ in clockwise order. Note that $v''Cv = eCv$. By Lemma \ref{two_edge_cut}, we may assume that $uCv'$ is contained in a block of $G - v''Cv$, which we denote by $H$. Since $|uCe| \ge 3$, $H$ is $2$-connected. See Figure \ref{fig:1}. Let $D$ be the outer cycle of $H$. By Lemma \ref{planar_technical_refined}, $H$ has a $D$-Tutte path $P_H$ between $u$ and $v'$, such that $f \in E(H)$ and

\begin{itemize}
	\item[(3)] $ \bridges[H]{P_H} \le (|H| - 6)/3 + \taubridges[H]{v'u} + \taubridges[H]{uf} + \taubridges[H]{fv'} - \bridgecorrection[H]{P_H}$.
\end{itemize}

Note that $\taubridges[H]{uf} = \taubridges[G]{uf}$ and $\taubridges[H]{fv'} = \taubridges[G]{fe}$. We now extend $P_H$ to the desired path in $G$. Let ${\cal B}$ be the set of $(H \cup v''Cv)$-bridges of $G$. Then $G = H \cup v''Cv \cup (\bigcup_{B \in {\cal B}}B)$. As $H$ is a block of $G - v''Cv$, we have $|B \cap H| \le 1$ for all $B \in {\cal B}$. Let ${\cal B}' = \{B \in {\cal B} : |B \cap H| = 1\}$ and define an equivalence relation on ${\cal B}'$: $B_1,B_2 \in {\cal B}'$ are equivalent if $B_1 \cap P_H = B_2 \cap P_H$, or if $B_1 \cap H$ and $B_2 \cap H$ are both disjoint from $P_H$ and  contained in the same $P_H$-bridge of $H$. Let ${\cal B}_i, i = 1, \ldots, m$, be the equivalence classes. We define $a_i,b_i \in V(v''Cv)$ such that
\begin{itemize}
	\item[(a)] $a_i \in V(B)$ and $b_i \in V(B')$ for some $B,B' \in {\cal B}_i$, allowing $B = B'$,
	\item[(b)] $v'', a_i,b_i,v$ occur on $eCv$ in order, and
	\item[(c)] subject to (a) and (b), $a_iCb_i$ is maximal.
\end{itemize}

Without loss of generality, we may assume that $a_1,b_1,a_2,b_2, \cdots, a_m,b_m$ occur on $v''Cv$ in order, with $v'' = a_1$ and $v = b_m$.  See Figure \ref{fig:1}. Let $J_i$ denote the union of $a_iCb_i$, all members of ${\cal B}_i$, all $(H \cup eCv)$-bridges of $G$ that have all attachments contained in $a_iCb_i$, and, if applicable, the $P_H$-bridge of $H$ containing $B \cap H$ for all $B \in {\cal B}_i$. Note that $|J_i\cap P_H|\in \{1,2\}$
for all $1\le i\le m$. Furthermore, for $1 \le i < m$, let $L_i$ denote the union of $b_iCa_{i+1}$ and the $(H \cup v''Cv)$-bridges of $G$ whose attachments are all contained in $b_iCa_{i+1}$. By (1), $|V(L_i)| \le 3$. Moreover, since $(G,C)$ is a circuit graph, $V(L_i) = V(b_iCa_{i + 1})$. Let

\begin{itemize}
	\item ${\cal J}_1$ be the set of all $J_i$ satisfying $|J_i\cap P_H|=1$  and  $a_i\ne b_i$, \sout{and}
	\item ${\cal J}_2$ be the set of all $J_i$ satisfying $|J_i\cap P_H|=2$.
\end{itemize}

By Lemma \ref{two_edge_cut}, $|J_1| = 2$ \track{and
	$v'' = a_1 = b_1 = a_2$}. Letting $P_1 = J_1$, we have

\begin{itemize}
	\item [(4)]  $\bridges[J_1]{P_1}=0=(|J_1|-1)/3 -1/3 - \bridgecorrection[J_1]{P_1}$.
\end{itemize}

Next, we consider ${\cal J}_1$ and ${\cal J}_2$. 

\begin{itemize}
	\item [(5)] For $J_i \in {\cal J}_1$,
	\sout{$J_i$}\track{$J_i - P_H$} has a path $P_i$ between $a_i$ and $b_i$, such that
	$P_i \cup (J_i\cap P_H)$ is an $a_iCb_i$-Tutte subgraph of $J_i$ and $\bridges[J_i]{P_i\cup (J_i\cap P_H)} \le (|J_i|-2)/3 - \bridgecorrection[J_i]{P_i\cup (J_i\cap P_H)}$. Furthermore, if $a_iCb_i$ is good, then 
\end{itemize}
\[
\bridges[J_i]{P_i\cup (J_i\cap P_H)}\le 
\begin{cases}
	(|J_i|-2)/3-1/3 -\bridgecorrection[J_i]{P_i\cup (J_i\cap P_H)},  & \text{if } |a_iCb_i| = 2,\\
	(|J_i|-2)/3-2/3 -\bridgecorrection[J_i]{P_i\cup (J_i\cap P_H)}, & \text{if } |a_iCb_i| > 2.
\end{cases}
\] 
To prove (5), let $V(J_i \cap P_H) = \{x\}$. Consider the plane graph $J_i' := J_i + a_ix$ whose outer cycle $C_i$ consists of $a_iCb_i$, the edge $e_i := xa_i$, and the path in the outer walk of $J_i$ between $b_i$ and $x$ not containing $a_i$. Note that $(J_i',C_i)$ is a circuit graph, $\taubridges[J_i']{xe_i} = 2/3$, and $\taubridges[J_i']{b_ix} = 0$. By Lemma \ref{planar_technical_refined}, $J_i'$ has a $C_i$-Tutte path $P_i'$ between $x$ and $b_i$, such that $e_i \in E(P_i')$ and $\bridges[J_i']{P_i'} \le (|J_i| - 6)/3 + \taubridges[J_i']{e_ib_i} + 2/3 - \bridgecorrection[J_i']{P_i'}$. Note that $\taubridges[J_i']{e_ib_i} \le 2/3$. If $a_iCb_i$ is good in $G$ then $\taubridges[J_i']{e_ib_i} \le 1/3$ (when $|a_iCb_i| = 2$) or $\taubridges[J_i']{e_ib_i} = 0$ (when $|a_iCb_i| > 2$). Hence, $P_i := P_i' - x$ gives the desired path for (5). \claimqed

\begin{itemize}
	\item [(6)] For $J_i \in {\cal J}_2$, $J_i - P_H$ has a path $P_i$ between $a_i$ and $b_i$, such that $P_i \cup (J_i \cap P_H)$ is an $a_iCb_i$-Tutte subgraph of $J_i$ and $	\bridges[J_i]{P_i\cup (J_i\cap P_H)} \le 	(|J_i|-4)/3 + 1 - \bridgecorrection[J_i]{P_i\cup (J_i\cap P_H)}$. Furthermore, if $|a_iCb_i| \ge 2$ and $a_iCb_i$ is good, then $	\bridges[J_i]{P_i\cup (J_i\cap P_H)}\le (|J_i|-4)/3+1/3 -\bridgecorrection[J_i]{P_i\cup (J_i\cap P_H)}$. 
\end{itemize}
It is easy to check that (6) holds trivially when $a_i = b_i$. So we may assume $a_i \neq b_i$. Let $V(J_i \cap P_H) = \{x,y\}$ such that $v',y,x,u$ occur on $D$ in clockwise order. Let $J_i'$ be the block of $J_i - \{x,y\}$ containing  $a_iCb_i$ and let $C_i$ be the outer cycle of $J_i'$. Note $a_iC_ib_i = a_iCb_i$. 

\track{If $J_i' \neq J_i - \{x,y\}$, we let $z \in V(J_i')$ denote the cut vertex \track{of $J_i - \{x,y\}$} separating $J_i' - z$ from $J_i - J_i' - \{x,y,z\}$. Otherwise,}
\sout{By planarity there exists a vertex}\track{by planarity, we let} $z \in V(b_iC_ia_i) - \{a_i,b_i\}$ such that $b_iC_iz - z$ contains no neighbor of $y$ and $zC_ia_i - z$ contains no neighbor of $x$. By Corollary \ref{planar_technical_vertex}, $J_i'$ has a $C_i$-Tutte path $P_i$ between $a_i$ and $b_i$, such that $z \in V(P_i)$ and 
\[ \bridges[J_i']{P_i} \le (|J_i'| - 3)/3 + \taubridges[J_i']{a_ib_i} - \bridgecorrection[J_i']{P_i}.\]

If $J_i' = J_i - \{x,y\}$ then $|J_i'| = |J_i| - 2$; so
\begin{align*}
	\bridges[J_i]{P_i \cup (J_i \cap P_H)} &= \bridges[J_i']{P_i}\\
	&\le (|J_i| - 2 -3)/3 + \taubridges[J_i']{a_ib_i} - \bridgecorrection[J_i']{P_i}\sout{.}\\
	&\sout{=}\track{\le} (|J_i| - 4)/3 + \taubridges[J_i']{a_ib_i} - 1/3 - \bridgecorrection[J_i]{P_i\cup (J_i\cap P_H)}\track{.}
\end{align*}

\sout{If $J_i' \neq J_i - \{x,y\}$ then $|J_i'| \le |J_i| - 3$; so}\track{If $J_i' \neq J_i - \{x,y\}$ then we have $|J_i'| \le |J_i| - 3$. Furthermore, as $(G,C)$ is a circuit graph, there are no nontrivial bridges of $P_i \cup (J_i \cap P_H)$ which either has $y$ and a vertex of $a_iC_iz - z$ as attachments or $x$ and a vertex of $zC_ib_i - z$ as attachments. Thus all vertices of $J_i - J_i' - \{x,y\}$ are contained in a single nontrivial bridge with $x,y$ and $z$ as attachments. It follows}
\begin{align*}
	\bridges[J_i]{P_i \cup (J_i \cap P_H)} &= \bridges[J_i']{P_i} + 1\\
	&\le (|J_i| - 3 -3)/3 + \taubridges[J_i']{a_ib_i} + 1 - \bridgecorrection[J_i']{P_i}\\
	&\sout{=}\track{\le} (|J_i| - 4)/3 + \taubridges[J_i']{a_ib_i} + 1/3 - \bridgecorrection[J_i]{P_i\cup (J_i\cap P_H)}.
\end{align*}

Now (6) holds, since $\taubridges[J_i']{a_ib_i} = 0$ if $a_iCb_i$ is good and $\taubridges[J_i']{a_ib_i} \le 2/3$ otherwise. \claimqed

Let $Q_i = b_iCa_{i+1}$ for $1 \le i < m$, and let $P = P_H \cup (\bigcup_{i=1}^m P_i) \cup (\bigcup_{i=1}^{m-1} Q_i)$. Then $P$ is a $C$-Tutte path in $G$ between $u$ and $v$ such that $\{e,f\} \subseteq E(P)$.  To bound $\bridges[G]{P}$, we count the $P$-bridges of $G$  in the following order: $\bridges[H]{P_H}$, $\bridges[J_1]{P_1}, \bridges[L_1]{Q_1}, \ldots, \bridges[L_{m-1}]{Q_{m-1}}, \bridges[J_m]{P_m}$.  Note that, for each $J_i \in {\cal J}_2$,  $|J_i\cap P_H|=2$ and $\bridges[G]{P}$ does not count the $P_H$-bridge of $H$ containing $J_i \cap H$. 
So we need to subtract $1$ from each $\bridges[J_i]{P_i \cup (P_H \cap J_i)}$ with $J_i \in {\cal J}_2$.
Hence, \begin{align*}
	\bridges[G]{P} &\le \bridges[H]{P_H} + \bridges[J_1]{P_1} + \sum_{J_i \in {\cal J}_1} \bridges[J_i]{P_i \cup (P_H \cap J_i)}\\
	&+ \sum_{J_i \in {\cal J}_2} (\bridges[J_i]{P_i \cup (P_H \cap J_i)} - 1) + \sout{\sum_{L_i \in {\cal L}}} \track{\sum_{i = 1}^{m-1}} \bridges[L_i]{Q_i}.
\end{align*}

Since $\bridges[L_i]{Q_i} = 0$ \track{for all $1 \le i \le m-1$}, we have

\begin{itemize}
	\item [(7)] $ \bridges[G]{P} \le \bridges[H]{P_H} + \bridges[J_1]{P_1} + \sum_{J_i \in {\cal J}_1 } \bridges[J_i]{P_i \cup (P_H \cap J_i)}$\\
	\hphantom{ssssssss}$ + \sum_{J_i \in {\cal J}_2} (\bridges[J_i]{P_i \cup (P_H \cap J_i)} - 1)$.
\end{itemize}

\sout{Note that for each $J_i\in {\cal J}_2$, $(|J_i\cap H|-3)/3$ is counted in $\bridgecorrection[H]{P_H}$ but not in $\bridgecorrection[G]{P}$. Hence,}  
\track{Furthermore, by (1), we have}
\begin{itemize}
	\item[(8)] $\sout{\bridgecorrection[G]{P}} \track{\bridgecorrection[G]{P} = 0} \le \sout{\bridgecorrection[H]{P_H} +} \bridgecorrection[J_1]{P_1} + \sum_{J_i \in {\cal J}_1} \bridgecorrection[J_i]{P_i \cup (P_H \cap J_i)} + \sum_{J_i \in {\cal J}_2} \bridgecorrection[J_i]{P_i \cup (P_H \cap J_i)}$.
\end{itemize}

We observe that

\begin{itemize}
	\item[(9)] $|J_1|$ and $|H|$ double count 1 vertex, namely $v'$, and, for $2 \le i \le m$ and $1 \le j \le m-1$, 
	\begin{itemize}
		\item[\textbullet]  $|L_j|$ and $|H \cup (J_1 \cup L_1) \cup \ldots \cup (J_{j-1} \cup L_{j-1}) \cup J_j|$ double count 1 vertex, namely \sout{$a_j$}\track{$b_j$},
		\item[\textbullet]  for $J_i \in {\cal J}_1$, $|J_i|$ and $|H \cup (J_1 \cup L_1) \cup \cdots \cup (J_{i-1} \cup L_{i-1})|$ double count 2 vertices, $a_i$ and the vertex in $V(J_i \cap P_H)$, and 
		\item[\textbullet]  for $J_i \in {\cal J}_2$, $|J_i|$ and $|H \cup (J_1 \cup L_1) \cup \cdots \cup (J_{i-1} \cup L_{i-1})|$ double count $|J_i \cap H| + 1$ vertices, $a_i$ and the vertices of  $J_i \cap H$.
	\end{itemize}
\end{itemize}

\track{If $J_i \in {\cal J}_2$, the vertices in $J_i \cap H - P_H$ are counted in both $|H|$ and $|J_i|$. Furthermore, as $\bridgecorrection[G]{P} = 0$, we have
	\[\bridgecorrection[H]{P_H} = \sum_{J_i \in {\cal J}_2} (|J_i \cap H| - 3)/3.  \] Thus for each $J_i\in {\cal J}_2$, $|J_i\cap H|-3$ is counted in $3\bridgecorrection[H]{P_H}$. Thus beginning with $H$ and $J_1$, and then iterating through $J_i$ and $L_i$ in the order they appear in $v''Cv$, by (9), we have the following,
	\begin{itemize}
		\item[(10)] \track{$|H| +  |J_1| - 1 + \sum_{J_i \in {\cal J}_1} (|J_i| - 2) + \sum_{J_i \in {\cal J}_2} (|J_i| - 4) + \sum_{i = 1}^{m-1} (|L_i| - 1) = n + 3\phi_H(P_H)$}.
	\end{itemize}
}
We now have three cases to analyze depending on the value of $\taubridges[G]{ev}$. In each case, we will find the desired path $P$.  Let ${\cal J}_2' = \{J_i \in {\cal J}_2 : a_i \neq b_i\}$, which, combined with (6), will be used in \sout{Case 2 and Case 3}\track{Cases 2 and 3}.\\

\textit{Case} 1. $\taubridges[G]{ev} = 2/3$.

If $eCv = v$, then  $P = P_H \cup v'Cv$ is a $C$-Tutte path in $G$ and $\{e,f\}\subseteq E(P)$. Moreover, $\bridgecorrection[G]{P}\sout{=}\track{\le}\bridgecorrection[H]{P_H}$; so by (3), 
\begin{align*}
	\bridges[G]{P} &\sout{=}\track{\le} \bridges[H]{P_H}\\
	&\le (|H| - 6)/e  + \taubridges[H]{v'u} + \taubridges[H]{uf} + \taubridges[H]{fv'} - \bridgecorrection[H]{P_H}\\
	&\le (|G| - 7)/e + 2/3 + \taubridges[G]{uf} + \taubridges[G]{fe} - \bridgecorrection[G]{P} \text{ (since $\taubridges[H]{v'u} \le 2/3$)}\\
	&= (n - 7)/3 + \taubridges[G]{ev} + \taubridges[G]{fe} + \taubridges[G]{uf} - \bridgecorrection[G]{P} \text{ (since $\taubridges[G]{ev} = 2/3$)}.
\end{align*}

If $eCv \neq v$ then $eCv$ is not good. Hence, by (3), (4), (5), (6), and (7), we have 
\begin{align*}
	\bridges[G]{P} &\le (|H| - 6)/3 + \taubridges[H]{v'u} + \taubridges[H]{uf} + \taubridges[H]{fv'} - \bridgecorrection[H]{P_H}\\
	&+ \sout{(}(|J_1| - 1)/3 - 1/3 -\bridgecorrection[J_1]{P_1}\sout{)}\\
	&+ \sum_{J_i \in {\cal J}_1}((|\sout{J_1}\track{J_i}| - 2)/3 -\bridgecorrection[J_i]{P_i \cup (P_H \cap J_i)})\\
	&+ \sum_{J_i \in {\cal J}_2} ((|J_i| - 4)/3 -\bridgecorrection[J_i]{P_i \cup (P_H \cap J_i)})   \\
	&\le (n - 7)/3 + \taubridges[G]{fe} + \taubridges[G]{uf} + \taubridges[H]{v'u} - \bridgecorrection[G]{P}  \text{ (by  \sout{(8) and (9)}\track{(8) and (10)})}\\
	&\le (n - 7)/3 + \taubridges[G]{fe} + \taubridges[G]{uf} + 2/3 - \bridgecorrection[G]{P} \text{ (since $\taubridges[H]{v'u} \le 2/3$)}\\
	&= (n - 7)/3 + \taubridges[G]{ev} + \taubridges[G]{fe} + \taubridges[G]{uf}  - \bridgecorrection[G]{P} \text{ (since $\taubridges[G]{ev} = 2/3$)}.\\
\end{align*}
\textit{Case} 2. $\taubridges[G]{ev} = 1/3$.

Then $|eCv| = 2$. Furthermore, $eCv$ is good. It follows from (3), (4), (5), (6), and (7) that 
\begin{align*}
	\bridges[G]{P} &\le (|H| - 6)/3  + \taubridges[H]{v'u} + \taubridges[H]{uf} + \taubridges[H]{fv'} - \bridgecorrection[H]{P_H}\\
	&+ \sout{(}(|J_1| - 1)/3 - 1/3 -\bridgecorrection[J_1]{P_1}\sout{)}\\
	&+ \sum_{J_i \in {\cal J}_1}((|J_i| - 2)/3 - 1/3 -\bridgecorrection[J_i]{P_i \cup (P_H \cap J_i)})\\
	&+ \sum_{J_i \in {\cal J}_2} ((|J_i| - 4)/3 -\bridgecorrection[J_i]{P_i \cup (P_H \cap J_i)}) -2|{\cal J}_2'|/3. 
\end{align*}
Thus, by \sout{(8) and (9)}\track{(8) and (10)}, 
\begin{align*}
	\bridges[G]{P} &\le  (n-7)/3 + \taubridges[G]{fe} + \taubridges[G]{uf} + \taubridges[H]{v'u}\\
	&- |{\cal J}_1|/3 - 2|{\cal J}_2'|/3 - |{\cal L}_2|/3 - 2|{\cal L}_3|/3 - \bridgecorrection[G]{P}\sout{,}\track{.}
\end{align*}
where ${\cal L}_2 = \{L_i : |L_i| = 2\}$ and ${\cal L}_3 = \{L_i : |L_i| = 3\}$.

As $|eCv| = 2$, $|{\cal J}_1| + |{\cal J}_2'| + |{\cal L}_2| + |{\cal L}_3|\ge 1$. Hence, since $\taubridges{ev} = 1/3$ and $\taubridges[H]{v'u} \le 2/3$, we have
$	\bridges[G]{P}  \le 	(n - 7)/3 + \taubridges[G]{uf}+ \taubridges[G]{fe}  + \taubridges[G]{ev} -\bridgecorrection[G]{P}.$\\

\textit{Case} 3. $\taubridges[G]{ev} = 0$.

First, suppose there exists $J_k \in {\cal J}_1$ such that $|a_kCb_k| > 2$. By (5), we may choose $P_k$ so that
\[ \bridges[J_k]{P_k \cup (J_k \cap P_H)} \le (|J_k| - 2)/3 - 2/3 - \bridgecorrection[J_k]{P_k \cup (J_k \cap P_H)}.\]
By (3), (4), (5), (6), and (7), we have
\begin{align*}
	\bridges[G]{P} &\le (|H| - 6)/3 + \taubridges[H]{v'u} + \taubridges[H]{uf} + \taubridges[H]{fv'} - \bridgecorrection[H]{P_H}\\
	&+ \sout{(}(|J_1| - 1)/3 - 1/3 -\bridgecorrection[J_1]{P_1}\sout{)}\\
	&+ \sout{(}(|J_k| -2)/3 -2/3 - \bridgecorrection{P_k \cup (J_k \cap P_H)}\sout{)}\\
	&+ \sum_{J_i \in {\cal J}_1 \setminus \{J_k\}}((|J_i| - 2)/3 - 1/3 -\bridgecorrection[J_i]{P_i \cup (J_i \cap P_H)})\\
	&+ \sum_{J_i \in {\cal J}_2} ((|J_i| - 4)/3 -\bridgecorrection[J_i]{P_i \cup (J_i \cap P_H)}) -2|{\cal J}_2'|/3  \\
	&\le (n-7)/3 + \taubridges[G]{uf} + \taubridges[G]{fe} + \taubridges[H]{v'u} - 2/3 -\bridgecorrection[G]{P} \text{ (by \sout{(8) and (9)}\track{(8) and (10)})}\\
	&\le (n-7)/3 + \taubridges[G]{uf} + \taubridges[G]{fe} + \taubridges[G]{ev} -\bridgecorrection[G]{P}.
\end{align*}

Now assume for all $J_k \in {\cal J}_1$, $|a_kCb_k| = 2$. Since $\taubridges[G]{ev} = 0$, $eCv$ is good and $|eCv| \ge 3$. Hence,
$|{\cal J}_1| + 2|{\cal J}_2'| + |{\cal L}_2| + 2|{\cal L}_3| \ge 2$.
By (3), (4), (5), (6), and (7), we have
\begin{align*}
	\bridges[G]{P} &\le (|H| - 6)/3 -1/3 + \taubridges[H]{v'u} + \taubridges[H]{uf} + \taubridges[H]{fv'} - \bridgecorrection[H]{P_H}\\
	&+ \sout{(}(|J_1| - 1)/3 - 1/3 -\bridgecorrection[J_1]{P_1}\sout{)}\\
	&+ \sum_{J_i \in {\cal J}_1;}((|J_i| - 2)/3 - 1/3 -\bridgecorrection[J_i]{P_i \cup (J_i \cap P_H)})\\
	&+ \sum_{J_i \in {\cal J}_2} ((|J_i| - 4)/3  -\bridgecorrection[J_i]{P_i \cup (J_i \cap P_H)}) -2|{\cal J}_2'|/3  \\
	&= (n-7)/3 + \taubridges[G]{uf}+ \taubridges[G]{fe}  + \taubridges[H]{v'u} \\
	&-|{\cal J}_1|/3 - 2|{\cal J}_2'|/3 - |{\cal L}_2|/3 - 2|{\cal L}_3|/3 -\bridgecorrection[G]{P}  \text{ (by  \sout{(8) and (9)}\track{ (8) and (10)})}\\
	&\le (n-7)/3 + \taubridges[G]{uf} + \taubridges[G]{fe} + \taubridges[G]{ev} -\bridgecorrection[G]{P}
\end{align*}
Hence, $P$ is the desired path. 

\section{Conclusion}\label{section_4}

\sout{For a positive integer $k$ and a graph $G$}\track{For a graph $G$ and a positive integer $k \ge 2$}, $G$ is \sout{said to be}\define{essentially} $k$-\define{connected} if $G$ is \track{$(k-1)$}-connected and, for every cut $S \subseteq V(G)$ of size less than $k$, $G - S$ has exactly two \sout{componenets}\track{components}, one of which is trivial. Quantitative results on Tutte paths lead to good bounds on length of longest cycles in essentially 4-connected graphs. For example, the authors \cite{WY20} showed that if $G$ is an essentially 4-connected planar graph on $n$ vertices, then it has a cycle of length at least $\ceil{(2n + 6)/3}$ and this bound is best possible. This result was independently proved by Kessler and Schmidt \cite{KS20} using a different method. 
As an immediate corollary of Theorem \ref{two_edge_tech}, we have the following

\begin{corollary}\label{cor:essentially_4_connected_cycle}
	Let $n \ge 4$ and let $G$ be an essentially 4-connected $n$-vertex plane graph. Let $e,f,g$ be distinct edges on \sout{a facial}\track{the outer} cycle of $G$. Then there is a cycle $C$ in $G$ such that $e,f,g \in E(C)$ and  $|C| \ge \ceil{ (2n + 1)/3}$. 
\end{corollary}
\track{\begin{proof}
		Let $u$ and $v$ be the endpoints of $g$ \sout{such that}\track{and assume} $u,f,e,v$ occur on $C$ in clockwise order. By Theorem \ref{two_edge_tech}, $G$ has a Tutte path $P$ between $u$ and $v$ such that $e,f \in E(G)$ and 
		\begin{align*}
			\bridges{P} &\le (n - 7)/3 + \taubridges{uf} + \taubridges{fe} + \taubridges{ev} - \bridgecorrection{P}\\
			&\le (n - 1)/3.
		\end{align*}
		Let $C = P + uv$. As $P$ is a Tutte path and $G$ is essentially 4-connected, we have that
		\begin{align*}
			|C| = n - \bridges{P} \ge (2n + 1)/3.
		\end{align*}
\end{proof}}

The bound on $|C|$ is tight with $K_4$ as an example. We conclude this section with \track{another} corollary of Theorem \ref{two_edge_tech}\track{.}\sout{, which will be used to obtain a quantitative result on Tutte paths in projective-planar graphs.} 

\begin{corollary}\label{two_edge_vertex}
	Let $n \ge 3$ be an integer, let $(G,C)$ be a circuit graph on $n$ vertices, let $e \in  E(C)$, and let $u,z,v \in V(C)$ be distinct, such that $u,z,e,v$ occur on $C$ in clockwise order\sout{ and $uv \in E(C)$}\track{, $z$ is not incident with $e$, and $uv \in E(C)$}. Then $G$ has a $C$-Tutte path $P$ between $u$ and $v$, such that $z \in V(P)$, $e \in E(P)$, and 
	\[\bridges{P} \le (n - 2)/3 - \bridgecorrection{P}.\]
\end{corollary}

\begin{proof}
	
	First, assume $G$ has no 2-cut contained in $uCe$. Then $\taubridges[G]{ue} = 0$ and any $C$-Tutte path in $G$ between $u$ and $v$ must also contain $z$. Applying Lemma \ref{planar_technical_refined}, we find a $C$-Tutte path $P$ in $G$ between $u$ and $v$, such that $e \in E(P)$ and
	\begin{align*}
		\bridges[G]{P} \le (n - 6)/3 + \taubridges[G]{vu} + \taubridges[G]{ue} + \taubridges[G]{ev} -\bridgecorrection[G]{P} \le (n - 2)/3  -\bridgecorrection[G]{P}.
	\end{align*}
	Since $z\in V(P)$, $P$ gives the desired path. 
	
	Now suppose there is a 2-separation $(G_1,G_2)$ in $G$ such that $\{u,v\} \subseteq V(G_1)$, $e \in E(G_1)$, and $V(G_1 \cap G_2) \subseteq V(uCe)$. We choose $(G_1,G_2)$ such that $z \in V(G_2) \setminus V(G_1)$ if possible. Let $V(G_1 \cap G_2) = \{x,y\}$ such that $u,x,y,e,v$ occur on $C$ in clockwise order. Let $G_1' = G_1 + xy$, and $C_1 := yCx + xy$ such that $(G_1',C_1)$ is a circuit graph. Let $e_1 := xy$. By Theorem \ref{two_edge_tech}, $G_1'$ has a $C_1$-Tutte path $P_1$ between $u$ and $v$, such that $\{e_1,e\} \subseteq E(P_1)$ and
	\begin{align*}
		\bridges[G_1']{P_1}  &\le (|G_1'| - 7)/3 + \taubridges[G_1']{ue_1}+ \taubridges[G_1']{e_1e} + \taubridges[G_1']{ev}  - \bridgecorrection[G_1']{P_1}\\
		&\le (|G_1| - 1)/3 - \bridgecorrection[G_1']{P_1}.
	\end{align*}
	Note that if $z \in V(G_1)$ then, by the choice of $(G_1,G_2)$, $z \in V(P_1)$. 
	Let $G_2' := G_2 + xy$ and $C_2 := xCy + yx$ such that $G_2'$ is a plane graph with outer cycle $C_2$. Then $(G_2',C_2)$ is a circuit graph. \sout{By Corollary \ref{planar_technical_vertex}, $G_2'$ has a $C_2$-Tutte path $P_2$ between $x$ and $y$ such that $z \in V(P_2)$ when $z \in V(G_2) \setminus V(G_1)$, and}\track{Let $z' = z$ if $z \in V(G_2) \setminus V(G_1)$, otherwise, let $z' \in V(C_2) \setminus \{x,y\}$. By Corollary \ref{planar_technical_vertex}, $G_2'$ has a $C_2$-Tutte path $P_2$ between $x$ and $y$ such that $z' \in V(P_2)$, and}
	\begin{align*}
		\bridges[G_2']{P_2} \le (|G_2'| - 3)/3 -  \bridgecorrection[G_2']{P_2} = (|G_2| - 3)/3 - \bridgecorrection[G_2']{P_2}.
	\end{align*}
	
	Now $P = (P_1 - xy) \cup P_2$ is a $C$-Tutte path in $G$ between $u$ and $v$, such that $z \in V(P)$ and $e \in E(P)$. Since $|G_1| + |G_2| = n + 2$  , we have
	\begin{align*}
		\bridges[G]{P} &\le (|G_1| - 1)/3 -\bridgecorrection[G_1']{P_1} + (|G_2| - 3)/3 - \bridgecorrection[G_2']{P_2}\\
		&= (n - 2)/3 - \bridgecorrection[G]{P}.
	\end{align*}
	So $P$ is the desired path. 
\end{proof}

\textbf{Acknowledgements} This work was completed while the first author was a graduate student at Georgia Institute of Technology. The authors would like to thank the anonymous referees for their careful reading and useful suggestions which improved the quality of the manuscript.

\end{document}